\documentclass[11pt, a4paper]{article}
\usepackage{a4wide} 
\usepackage{graphicx} % Required for inserting images
\usepackage[T1]{fontenc}
\usepackage[utf8]{inputenc}
\usepackage{lmodern} 
\usepackage{amssymb,amsmath,bm, nccmath}
\usepackage{graphicx}
\usepackage{tikz}
\usetikzlibrary{decorations.pathreplacing}
\usepackage{pgfplots}
\usepackage{subcaption}
\usepackage{multicol}

\usepackage{booktabs}

\newcommand{\R}{\mathbb{R}}

\newcommand{\E}{\mathcal{E}}
\newcommand{\V}{\mathcal{V}}
\newcommand{\J}{\mathcal{J}}
\newcommand{\Scal}{\mathcal{S}}
\newcommand{\Vmax}{V^{\max}}
\newcommand{\Vbm}{\bm{V}^{\max}}
\newcommand{\Vfix}{\bm{V}^{\mathrm{fix}}}
\newcommand{\x}{\bm{x}}
\newcommand{\vv}{\bm{v}}
\newcommand{\alphaa}{\bm{\alpha}}
\newcommand{\alphafix}{\bm{\alpha}^{\mathrm{fix}}}

\newcommand{\dd}{\, \mathrm{d}}

\newcommand{\inn}{\mathrm{in}}
\newcommand{\out}{\mathrm{out}}
\newcommand{\flow}{\mathrm{flow}}
\newcommand{\diff}{\mathrm{diff}}

\newcommand{\opt}{\mathrm{opt}}
\newcommand{\road}{{\mathrm{road}}}

\newcommand{\red}[1]{\textcolor{red}{#1}}

\begin{document}

\title{Influence of Routing and Speed Limits on Optimal Solutions in Traffic Emission Modeling}
\author{Marc-André Bach\footnotemark[1], \; Simone G\"ottlich\footnotemark[1], \;Alena Ulke\footnotemark[1]} 

\footnotetext[1]{University of Mannheim, Department of Mathematics, 68131 Mannheim, Germany (goettlich@uni-mannheim.de, ulke@uni-mannheim.de)}

\maketitle   

\begin{abstract}
We investigate the influence of routing strategies and speed limit policies on optimal solutions in traffic emission models. Building on a first-order macroscopic traffic model coupled with an advection–diffusion model, we formulate single- and multi-objective optimization problems to simultaneously maximize traffic efficiency and minimize air pollution. 
We compare three control scenarios: optimizing only the routing strategy, optimizing only the speed limit policy, and optimizing both simultaneously. Numerical experiments on a small road network demonstrate that speed limit policies consistently achieve larger reductions in emissions and greater gains in traffic efficiency than routing strategies. Multi-objective optimization reveals the trade-off between the two goals and confirms that including speed limits in the control set yields Pareto-optimal solutions that are strictly superior to those obtained by routing control only. Our results provide quantitative guidance for traffic management seeking to balance mobility and environmental objectives.
\end{abstract}

\section{Introduction and Motivation}

The reduction of greenhouse gas and air pollutant emissions from road traffic is one of the central challenges of urban mobility planning. In many metropolitan areas, vehicular traffic contributes substantially to local concentrations of nitrogen oxides, particulate matter, and other harmful substances, making traffic management a key lever for improving air quality \cite{WHO2021}. 
At the same time, efficient traffic flow is an important economic and social objective: congestion leads to longer travel times, higher fuel consumption, and hence increased emissions.
Two classes of control instruments that are available to manage vehicular traffic are speed limit policies, which set the maximum allowed velocity on individual road segments, and routing strategies, which influence the distribution of traffic across the network by prescribing the use of particular roads. Both affect emissions through their impact on traffic density and flow, but their effectiveness has not been compared yet.

%\medbreak
Macroscopic traffic models such as the Lighthill-Whitham-Richards (LWR) model~\cite{LighthilWhitham1955,Richards1956} or the Aw-Rascle-Zhang (ARZ) model~\cite{AwRascle2000,Zhang2002} provide a basis for the analysis and optimization of large-scale traffic networks.
Building on these models, traffic dynamics have been extended to networks that include junctions and require routing strategies and priority rules at intersections, cf. \cite{Bressan2014,Garavello2016,Piccoli2006}. 
The environmental impact of traffic, in particular air pollution, has also been modeled and investigated in \cite{Alvarez2017,Balzotti2022,Berrone2012}, where the authors couple a traffic model with an advection-diffusion model to describe the evolution of the concentration of pollutants in the atmosphere.
These approaches have been further extended in \cite{Alvarez2018,Alvarez2019} to optimization problems that aim to reduce the environmental impact of traffic.
Also multi-objective optimization of traffic networks has attracted growing interest. For example,\cite{Goatin2026} addresses the trade-off between minimizing the total travel time and the total travel distance in a multi-class traffic model while controlling the routing strategy in the network.
Furthermore, \cite{Ulke2025} adopts a multi-objective framework to analyze the trade-off between emission reduction and increasing traffic efficiency by controlling speed limits. 
We build on the framework of \cite{Ulke2025} and additionally emphasize the influence of routing parameters at junctions within the traffic model. Moreover, we extend their analysis by comparing the effect of three control strategies on the 
on the optimal solutions of single- and multi-objective optimization problems: controlling only the routing strategy, only the speed limit policy, and controlling both simultaneously.
Our goal is to provide insights for decision makers on whether routing or speed limit control offers greater potential to reduce emissions and increase traffic efficiency.

The remainder of this paper is organized as follows. Section 2 recalls the traffic and emission model. In Section 3, we formulate suitable optimization problems, which we solve numerically in Section 4 for a test network and analyze the results. %Section 5 draws conclusions and outlines future research directions.

\section{Modeling Emissions of Vehicular Traffic}\label{sec: traffic emission model}

Modeling the impact of traffic on both the environment and traffic efficiency involves three components:
(i) the traffic dynamics on a road network, including the influence of routing strategies and speed limit policies, 
(ii) the estimation of air pollutant emissions generated by traffic, and
(iii) the transport and dispersion of these pollutants in the atmosphere due to wind and diffusion.
Thus, the overall framework consists of three coupled models: a traffic model, an emission model, and a dispersion model.
We adopt the modeling framework proposed in \cite{Alvarez2017,Alvarez2018}, which is extended to incorporate speed limit policies in \cite{Ulke2025}. In the following, we briefly introduce the individual components and highlight where the routing strategy enters the model.

\subsection{Traffic Model}

A traffic network is a directed graph \(\mathcal{G} = (\V, \E)\) where each edge \(e \in \E\) models a unidirectional road and is associated with an interval \( I_e\), and each node \(v \in \V\) models a junction connecting roads.
We model the traffic flow on an edge \(e\) using the \emph{Lighthill-Whitham-Richards (LWR) model}~\cite{LighthilWhitham1955,Richards1956} which reads
\begin{subequations}\label{LWRmodel}
	\begin{align}
		\partial_t\rho_e(s,t)+\partial_s Q_e(\rho_e(s,t),\Vmax_e) &= 0, &&\text{for } (s,t) \in I_e \times (0,T),\\
		\rho_e(s,0) &= \rho_e^0(s), &&\text{for } s \in I_e,
	\end{align}
\end{subequations}
where \(\rho_e\) denotes the traffic density, \(\rho_e^{\max}\) the maximal traffic density or the road's capacity, and  \(\Vmax_e\) the speed limit. 
Further, \(Q_e\) describes the flux function, for which we assume Greenshields' model~\cite{Greenshields1935}, i.e.,
\begin{align*}
	Q_e(\rho, \Vmax_e) & =  \Vmax_e \rho \left( 1 - \frac{\rho}{\rho^{\max}_e} \right).
\end{align*}
The LWR model admits a weak solution on every road if we fix the speed limits and provide additional boundary data, cf.~\cite{Piccoli2006,Garavello2016}. As traffic on different roads interacts at junctions, the solution depends on the speed limit policy \(\Vbm = (\Vmax_e)_{e\in\E}\) of the entire network. Thus, we write \(\rho_e(s,t) = \rho_e(s,t,\Vbm)\) in the following.

%\medbreak
Traffic flow across junctions is governed by imposing \emph{suitable} boundary conditions at the beginning and end of each road. These boundary conditions either arise from coupling conditions or from prescribing boundary data explicitly.
We distinguish between external junctions where drivers enter or exit the network, and internal ones connecting multiple roads.% \red{cf. Figure \ref{fig: network} for an illustration.}

At an external junction where drivers enter the network, we provide an inflow rate \(q_e^{\inn}(t)\) of vehicles entering the road \(e\), cf.~\cite{Goatin2016,Herty2009}, and at an external junction where drivers exit the network, we assume free flow.

An internal junction can be seen as a tuple \((e_1,\ldots,e_n,e_{n+1}, e_{n+m})\) where the first \(n\) entries correspond to incoming roads and the last \(m\) entries to outgoing ones.
At each internal junction, we impose mass conservation together with additional rules to characterize how traffic is distributed from incoming to outgoing roads.
We follow \cite{Piccoli2006,Garavello2016} and assume (i) that traffic is distributed according to drivers preferences, which are modeled by a \emph{preference or routing parameter}, (ii) at merging junctions, where \(n < m\) holds, the available capacity of the outgoing roads is allocated based on the priority of the incoming roads, which is encoded by a \emph{priority parameter}, and (iii) drivers aim to maximize the flow across the junction respecting these rules.

In order to formalize these coupling conditions, we require the demand and supply function for each road, respectively:
\begin{align*}
	D_e(\rho, \Vmax_e) = \begin{cases}
		Q_e(\rho, \Vmax_e), & \rho \le \rho_e^c, \\
		Q_e^{\max}(\Vmax_e), & \rho > \rho_e^c,
	\end{cases} \quad \text{and} \quad 
	S_e(\rho, \Vmax_e) = \begin{cases}
		Q_e^{\max}(\Vmax_e), & \rho \le \rho_e^c, \\
		Q_e(\rho, \Vmax_e), & \rho > \rho_e^c,
	\end{cases}
\end{align*}
where \(\rho_e^c\) is the critical traffic density for which the flux function \(Q_e\) reaches its maximum \(Q_e^{\max}\). As we assume Greenshields' model for the flux function, we have \(\rho_e^c = \rho_e^{\max}/2\). Furthermore, we introduce the following abbreviations:
\begin{align*}
	D_i^{\inn}(\rho, \Vbm) & = D_{e_i}(\rho_{e_i}(L_{e_i}^\road,t,\Vbm), \Vmax_{e_i})
		&& \text{for} \quad i = 1, \ldots, n, \\
	S_j^{\out}(\rho, \Vbm) & = D_{e_j}(\rho_{e_j}(0,t,\Vbm), \Vmax_{e_j})
		&& \text{for} \quad j = n+1, \ldots, n+m,
\end{align*} 
where \(L_e^{\road}\) denotes the length of the road \(e\). Then the boundary conditions at an internal junction at a fixed time \(t \in [0,T]\) are given by the solution to the constrained maximization problem:
\begin{subequations}\label{maxProbJunc}
	\begin{align}
		\max_{Q_{e_i}^{\inn}(t,\Vbm)} 
			\,\,\, & \sum_{i=1}^n Q_{e_i}^{\inn}(t,\Vbm),\\
		\text{subject to} \quad  & \sum_{i=1}^n \alpha_{j,i} Q_{e_i}^{\inn}(t,\Vbm) \leq S_j^{\out}(t,\Vbm), 
			&& \text{for all} \quad j = n+1,\ldots,n+m, \label{supCond}\\
		& 0 \leq Q_{e_i}^{\inn}(t,\Vbm) \leq D_i^{\inn}(t,\Vbm), 
			&& \text{for all} \quad i = 1,\ldots,n.\label{demCond}
	\end{align}
\end{subequations}
The parameters \(\alpha_{j,i} \in [0,1]\) are the so-called \emph{preference or routing parameters} and encode the drivers' preference to take the outgoing road \(j\) from the incoming road \(i\) and must satisfy \( \sum_{i=n+1}^{n+m} \alpha_{j,i} = 1\).
The explicit solution to \eqref{maxProbJunc} in case of a one-to-one, a merging two-to-one, and a diverging one-to-two junction can be found in \cite{Goatin2016}. For general \(n\)-to-\(m\) junctions the existence of a solution is discussed in \cite{Piccoli2006,Garavello2016} and the references therein.
Since the traffic dynamic on the roads is affected by the dynamic at the junctions, the traffic densities also depend on the routing parameter or routing strategy \(\alphaa = \left(\alpha_{j,i}\right)_{j,i} \in \R^{m\times n}\). Hence, we expand our notation to \(\rho_e(s,t) = \rho_e(s,t,\alphaa, \Vbm)\).

%At external junctions, we provide boundary data in form of incoming or outgoing flow. However, as the imposed flow might exceed the roads available capacity, we set:
%\begin{equation}\label{inflow}
%	Q_{e}^{in} = \min\lbrace S_e^{in}(t,\boldsymbol{V}^{max}),q_e^{in}(t) \rbrace.
%\end{equation}
%\begin{equation}\label{eqOutflow}
%	Q_{e}^{out} = \min\lbrace D_e^{in}(t,\boldsymbol{V}^{max}),f_{e}^{out}(t) \rbrace.
%\end{equation}
%where \(D_e\) and \(S_e\) are the supply and demand function.

\subsection{Emission and Dispersion Model}

After providing a model to describe the traffic dynamic on a network, we continue by modeling the contribution of traffic to air pollution as follows:
(i) we model the emission of pollutants by traffic itself and then (ii) we model the spread of pollutants in a given area in the air where we will utilize the emission estimate from (i).

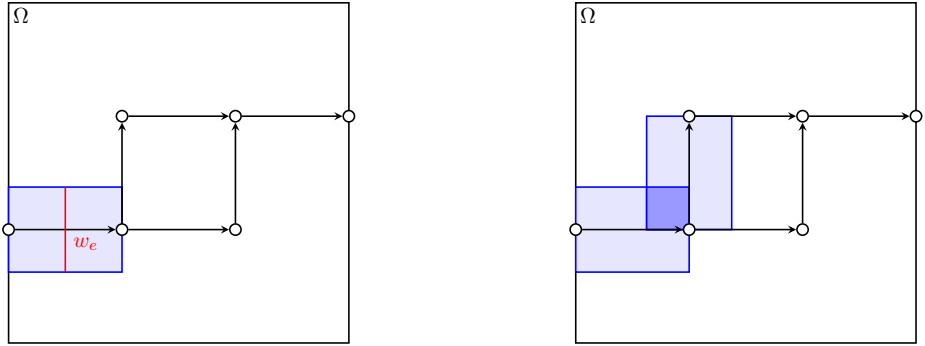
\begin{figure}[h]
	\centering
    \scalebox{0.75}{
	\begin{tikzpicture}[node distance=2cm, >=stealth, auto,thick]
        \draw[thick] (0,-2) rectangle (6,4);
        \node[anchor=south east] at (0.5,3.5) {$\Omega$};
    
        \draw[thick, blue, fill=blue!10] (0,-0.75) rectangle (2,0.75);
                    
        \node (1) [circle, draw=black, fill=white, inner sep=2pt] {};
        \node (2) [circle, draw=black, fill=white, inner sep=2pt, right of=1] {};
        \node (3) [circle, draw=black, fill=white, inner sep=2pt, above of=2] {};
        \node (4) [circle, draw=black, fill=white, inner sep=2pt, right of=2] {};
        \node (5) [circle, draw=black, fill=white, inner sep=2pt, above of=4] {};
        \node (6) [circle, draw=black, fill=white, inner sep=2pt, right of=5] {};

        % Draw arrows with labels
        \draw[->] (1) -- (2) node[midway, below] {};
        \draw[->] (2) -- (3) node[midway, left] {};
        \draw[->] (2) -- (4) node[midway, below] {};
        \draw[->] (3) -- (5) node[midway, above] {};
        \draw[->] (4) -- (5) node[midway, right] {};
        \draw[->] (5) -- (6) node[midway, above] {};

        \draw[thick, red] (1,-0.75) -- (1,0.75) node[midway, below right, red] {$w_e$};
        
         \begin{scope}[xshift=10cm]
            \draw[thick] (0,-2) rectangle (6,4);
            \node[anchor=south east] at (0.5,3.5) {$\Omega$};
    
            \draw[thick, blue, fill=blue!10] (0,-0.75) rectangle (2,0.75);
            \draw[thick, blue, fill=blue!10] (1.25,0) rectangle (2.75,2);
                                    
            \draw[thick, blue, fill=blue!40] (1.25,0) rectangle (2,0.75);

            \node (1) [circle, draw=black, fill=white, inner sep=2pt] {};
            \node (2) [circle, draw=black, fill=white, inner sep=2pt, right of=1] {};
            \node (3) [circle, draw=black, fill=white, inner sep=2pt, above of=2] {};
            \node (4) [circle, draw=black, fill=white, inner sep=2pt, right of=2] {};
            \node (5) [circle, draw=black, fill=white, inner sep=2pt, above of=4] {};
            \node (6) [circle, draw=black, fill=white, inner sep=2pt, right of=5] {};
    
            % Draw arrows with labels
            \draw[->] (1) -- (2) node[midway, below] {};
            \draw[->] (2) -- (3) node[midway, left] {};
            \draw[->] (2) -- (4) node[midway, below] {};
            \draw[->] (3) -- (5) node[midway, above] {};
            \draw[->] (4) -- (5) node[midway, right] {};
            \draw[->] (5) -- (6) node[midway, above] {};
         \end{scope}
    \end{tikzpicture}
    }
	\caption{The spatial layout of a road network in a two-dimensional domain \(\Omega\) including the extension of a road to its width \(w_e\) in blue (left) and overlapping of two extended roads (right). }
	\label{fig: emission rate}
\end{figure}

%\medbreak
Since we model the traffic dynamic using a first-order model, we follow \cite{Alvarez2017,Alvarez2018} and estimate emissions by a linear combination of traffic density and flow. 
For a given weighting parameter \(\theta\), the emission rate \(\xi_e\) on road \(e\) is defined by
\begin{equation}\label{eq: emission rate}
    \xi_e(s, t, \alphaa,\Vbm) = Q_e(\rho_e(s, t, \alphaa,\Vbm),\Vmax_e) + \theta\rho_e(s, t, \alphaa, \Vbm),
\end{equation}
where \((s, t) \in I_e \times (0,T)\). As we aim to model the spread of pollutants in a two-dimensional area \(\Omega\), we need to transfer the emission rates \(\xi_e\) from the one-dimensional road network to \(\Omega\).
We follow the approach presented in \cite{Ulke2025} and map each road \(e\) to its spatial location in \(\Omega\) and extend it to a finite width \(w_e\). 
Then, we distribute the emissions uniformly across the roads' width. 
However, the extended roads may overlap. In such regions, we define the emission rate by the average of the emission rates of overlapping roads. 
This results in the emission rate \(\xi \colon \Omega \times [0,T] \rightarrow \R\) with \((\x,t) \mapsto \xi(\x,t,\alphaa,\Vbm)\), which defined on the two-dimensional area \( \Omega\), see Figure~\ref{fig: emission rate} for an illustration. 
%\begin{equation*}
	%\xi \colon \Omega \times [0,T] \rightarrow \R, \quad (\x,t) \mapsto \xi(\x,t,\alphaa,\Vbm).
%\end{equation*}

%\medbreak
In addition to the emission of air pollutants due to traffic, we are also interested in the evolution of the pollutant concentration \(\phi\). As the driving factors of the dispersion are wind and diffusion, we employ an advection-diffusion equation to model this phenomenon.
In particular, we assume that traffic is the only source of pollution within \(\Omega\). Moreover, pollutants may leave the domain due to the wind dynamic, but cannot enter it from outside. To formalize this assumption, we introduce the inflow and outflow boundary of \(\Omega\) for a possibly space- and time-dependent wind field \(\vv = (v_x,v_y)^\top\), denoted by \(S^-(\vv)\) and \( S^+(\vv)\), respectively:
\begin{align*}
	S^- (\vv) 
		& = \lbrace (\x,t) \in \partial\Omega \times [0,T] \mid \vv(\x,t) \cdot \overset{\rightarrow}{\eta}(\x) < 0\rbrace, \\
	S^+ (\vv) 
		& = \lbrace (\x,t) \in \partial\Omega \times [0,T] \mid \vv(\x,t) \cdot \overset{\rightarrow}{\eta}(\x) \geq 0\rbrace,
\end{align*}
where $\overset{\rightarrow}{\eta}$ is the unit outward normal of \(\Omega\).
For a given initial concentration \(\phi_0\) of pollutants, the evolution of \(\phi\) is then governed by the following advection-diffusion equation, cf. \cite{Alvarez2017,Balzotti2022,Stockie2011}:
\begin{subequations}\label{eq: IBVP}
	\begin{align}
		\partial_t\phi(\x,t) - \mu\Delta\phi(\x,t) + \vv\cdot\nabla\phi(\x,t) + \kappa \phi(\x,t) & =  \xi(\x,t,\alphaa, \Vbm), 	
			&&(\x,t) \in \Omega \times (0,T), \label{eq: advection diffusion}\\
		\mu\nabla\phi(\x,t) \cdot\overset{\rightarrow}{\eta}(\x) - \phi(\x,t)\vv\cdot\overset{\rightarrow}{\eta}(\x) & = 0, 
			&&(\x,t) \in S^-(\vv), \label{eq: inflow}\\
		\nabla\phi(\x,t) \cdot\overset{\rightarrow}{\eta}(\x) & = 0, && (\x,t) \in S^+(\vv), \label{eq: outflow}\\
		\phi(\x,0) & = \phi_0(\x), && \x \in \Omega.
	\end{align}
\end{subequations}
The constant \(\mu\) describes the diffusion coefficient, the constant \(\kappa\) the extinction rate of the air pollutants, and the source term in Equation \eqref{eq: advection diffusion} accounts for emissions generated by traffic. 
Furthermore, the boundary conditions \eqref{eq: inflow} and \eqref{eq: outflow} ensure that no pollutants enter \(\Omega\) from the outside, but may leave \(\Omega\) due to wind and diffusion.
The existence of a weak solution to the initial boundary value problem~\eqref{eq: IBVP} was established in \cite{Skiba2000}.

\section{Optimization}\label{sec: optimization}

Our main goal is to assess the impact of routing strategies and speed limit policies on traffic efficiency and traffic emissions. In particular, our goal is to minimize emissions and maximize efficiency while investigating the effectiveness of controlling routing, speed limits, or both on optimal solutions.
Therefore, we introduce suitable cost functions that quantify the environmental impact and efficiency of traffic. 
Based on these cost functions, we formulate optimization problems that either minimize the environmental impact or maximize traffic efficiency. Additionally, we consider a multi-objective optimization framework where both objectives are optimized simultaneously, allowing us to analyze trade-offs between the objectives.

\medbreak
We measure the efficiency of traffic using the so-called \emph{accumulated traffic flow}, which aggregates the traffic flow on all roads of the network over the entire time horizon, cf.~\cite{Goettlich2015,Ulke2025,Treiber2013}:
\begin{equation*}\label{Jflow}
	\J_{\flow}( \alphaa, \Vbm) :=  -\sum_{e \in \E}\int_0^T\int_{I_e} Q_e(\rho_e(s,t,\alphaa,\Vbm),\Vmax_e) \dd s \dd t.
\end{equation*}
We introduce the negative sign because the optimization will be formulated as a minimization problem, whereas we want to maximize the traffic efficiency.
To quantify the environmental impact of traffic emissions, we employ the \emph{mean contamination}, which measures the average pollutant concentration over the spatial domain and the considered time horizon, cf.~\cite{Alvarez2017,Alvarez2018,Ulke2025}:
\begin{equation*}\label{Jdiff}
	\J_{\diff}(\alphaa,\Vbm) := \frac{1}{T \vert\Omega\vert}\int_0^T\int_{\Omega} \phi(\x,t,\alphaa,\Vbm) \dd \x \dd t.
\end{equation*}
Evaluating mean contamination for different routing strategies and speed limit policies requires solving a coupled system of partial differential equations (PDEs) repeatedly.
Since algorithms in numerical optimization typically require many evaluations of the cost function, this become computationally costly. 
To reduce the number of PDE solves, we employ adjoint calculus, which allows to express the cost function \( \J_\diff\) in terms of an adjoint variable \(p\), the source term \(\xi\), and the initial data \(\phi_0\):
\begin{equation*}%\label{Jdiff}
	\J_{\diff}(\alphaa,\Vbm) := \int_0^T \int_{\Omega} \xi(\x,t,\alphaa, \Vbm) p(\x,t) \dd \x \dd t
		+ \int_{\Omega} \phi_0(\x) p(\x,0) \dd \x.
\end{equation*}
The adjoint \(p\) is obtained by solving an advection-diffusion equation backwards in time, cf. \cite{Alvarez2018} for details. A key advantage of this approach is that the adjoint equation is independent of the routing parameter \(\alphaa\) and the speed limits \(\Vbm\).
Thus, once the adjoint \(p\) is computed, evaluating \(\J_\diff\) for different control parameters only requires recomputing the source term \(\xi\), which in turn depends only on the solution to the traffic model.
Hence, we eliminate the need to solve the dispersion model repeatedly.

\medbreak
After introducing the cost functions, our next step is to formulate corresponding optimization problems.
We first consider the standard single-objective setting before extending the framework to multi-objective optimization.

Depending on the choice of the control variables, we need to impose suitable constraints.
In the case where only the routing strategy \(\alphaa\) is controlled, we prescribe lower and upper bounds \(\underline{\alpha}_{j,i} \) and \( \overline{\alpha}_{j,i}\) for the routing parameters.
Then, the standard optimization problem takes the form:
\begin{equation*}\label{singleOptalpha}
	\min_{\alphaa \in \Scal_{\alphaa}} \J_{\ast}(\alphaa, \Vfix) \quad \text{where} \quad 
		\Scal_{\alphaa} =  \Bigl\lbrace \alphaa \in [0,1]^{k} \mid \underline{\alpha}_{j,i} \leq \alpha_{j,i} \leq \overline{\alpha}_{j,i}\Bigr\rbrace
\end{equation*}
where \(k\) is the number of routing parameters and \(\ast \in \{\flow, \diff\}\).
As we only control the routing parameter, the speed limit policy need to be fixed explicitly. We denote this fixed speed limit policy by \(\Vfix\).
In case where we only control the speed limit policy \(\Vbm\) but not the routing strategy \(\alphaa\), we also impose lower and upper bounds \(\underline{V_e}\) and \( \overline{V}_e\), respectively. Then, the standard optimization problem reads:
\begin{equation*}\label{singleOptVmax}
	\min_{\Vbm \in \Scal_{\bm{V}}} \J_{\ast}(\alphafix,\Vbm) \quad \text{where} \quad 
	\Scal_{\bm{V}} =  \Bigl\lbrace \Vbm \in [0,\infty)^{d} \mid \underline{V_e} \leq \Vmax_e \leq \overline{V}_e\Bigr\rbrace
\end{equation*}
where \(d\) is the number of roads in the network and and \(\ast \in \{\flow, \diff\}\).
Again, as we only control the speed limit policy, we need to provide a fixed routing strategy \(\alphafix\) a priori.
Finally, combing both approaches and controlling the routing strategy and the speed limit policy simultaneously, yields the following:
\begin{subequations}\label{singleOptBoth}
	\begin{align*}
		\min_{(\alphaa, \Vbm) \in \Scal} &  \,\, \J_{\ast}(\alphaa, \Vbm) \\
		\text{where} \quad & \Scal = \Bigl\lbrace (\alphaa,\Vbm) \in [0,1]^{k}\times[0,\infty)^{d} \mid 
			\underline{\alpha}_{j,i} \leq \alpha_{j,i} \leq \overline{\alpha}_{j,i}, \,\, 
			\underline{V_e} \leq \Vmax_e \leq \overline{V}_e\Bigr\rbrace,
	\end{align*}
\end{subequations}

\medbreak
Since both emission reduction and traffic efficiency are important objectives in traffic management, we also consider optimizing them simultaneously to investigate the trade-off. 
Therefore, we formulate the following multi-objective optimization problem:
\begin{equation}\label{eq: multiobjective opt}
	\min_{(\Vbm,\alphaa) \in \Scal} \,\,\, \bm{\J}(\bm{V}) = \left(\J_\flow, \J_\diff\right).
\end{equation}
Analogously, with respect to the standard optimization problems with one single objective, optimization may be considered with respect to only the routing strategy or only the speed limit policy. In this case, we replace the feasible set \(\Scal\) by  \(\Scal_{\alphaa}\) or \(\Scal_{\bm{V}}\) respectively, while the remaining control is fixed.

In multi-objective optimization, the individual objectives generally do not attain their minimum for the same control. Such objectives are also called \emph{conflicting}, meaning that there exists no control that minimizes all objectives simultaneously.
Hence, we characterize optimality using the concept of Pareto optimality:
a control \((\Vbm_\opt, \alphaa_{\opt})\) is called \emph{Pareto optimal} if no objective value can be improved without deteriorating at least one of the others. The corresponding objective value \(\bm{\J}(\Vbm_\opt, \alphaa_{\opt})\) is called \emph{efficient} and lies on the Pareto front. For comprehensive introductions to multi-objective optimization and Pareto optimality, we refer to \cite{Branke2008,Ehrgott2005}.
Solving the multi-objective optimization problem \eqref{eq: multiobjective opt} means determining the Pareto front.

\section{Numerical Examples}

After formulating the optimization problems to minimize pollution and maximize traffic efficiency, we now investigate their solutions for a test case. Our goal is to assess the influence of routing strategies and speed limit policies on both traffic efficiency and environmental impact. In particular, our objective is to determine which control mechanism has the greatest effect and to compare the resulting optimal routing strategies and speed limit policies.

\medbreak
Since the PDEs that govern the traffic emission model generally lack closed-form solutions, we solve the optimization problems numerically using a \emph{first-discretize-then-optimize} approach.
We employ the discretization proposed in \cite{Ulke2025}: the traffic model is discretized using Godunov's scheme, while the adjoint advection--diffusion equation is discretized by finite differences.
Finally, we use the right-rectangular rule to resolve the integrals in the objectives. The resulting discrete objective functionals are denoted by \(\J_{\flow}^h\) and \(\J_{\diff}^h\).

\begin{figure}[htb!]
	\begin{center}
		\scalebox{0.75}{
		\begin{tikzpicture}[node distance=2cm, >=stealth, auto, thick]
			
			\draw[thick] (0,-2) rectangle (6,4);
			\node[anchor=south east] at (0.5,3.5) {$\Omega$};

			\node (1) [circle, draw=black, fill=white, inner sep=2pt] {};
			\node (2) [circle, draw=black, fill=black, inner sep=2pt, right of=1] {};
			\node (3) [circle, draw=black, fill=white, inner sep=2pt, above of=2] {};
			\node (4) [circle, draw=black, fill=white, inner sep=2pt, right of=2] {};
			\node (5) [circle, draw=black, fill=white, inner sep=2pt, above of=4] {};
			\node (6) [circle, draw=black, fill=white, inner sep=2pt, right of=5] {};
			
			% Draw arrows with labels
			\draw[->] (1) -- (2) node[midway, below, inner sep=2pt] {$ e=1$};
			\draw[->] (2) -- (3) node[midway, right, inner sep=2pt] {$\red{\alphaa}$};
            \draw[->] (2) -- (3) node[midway, above,rotate=90, inner sep=2pt] {$e=2$};
			\draw[->] (2) -- (4) node[midway, above, inner sep=2pt] {$\red{\bm{1-}\alphaa}$};
            \draw[->] (2) -- (4) node[midway, below, inner sep=2pt] {$e=3$};
			\draw[->] (3) -- (5) node[midway, above, inner sep=2pt] {$e=4$};
			\draw[->] (4) -- (5) node[midway, below, rotate=90, inner sep=2pt] {$e=5$};
			\draw[->] (5) -- (6) node[midway, above, inner sep=2pt] {$e=6$};

            \begin{scope}[xshift=10cm]			
			    \draw[thick] (0,-2) rectangle (6,4);
			    \node[anchor=south east] at (0.5,3.5) {$\Omega$};

			    \node (1) [circle, draw=black, fill=white, inner sep=2pt] {};
			    \node (2) [circle, draw=black, fill=white, inner sep=2pt, right of=1] {};
		          \node (3) [circle, draw=black, fill=white, inner sep=2pt, above of=2] {};
		          \node (4) [circle, draw=black, fill=white, inner sep=2pt, right of=2] {};
		          \node (5) [circle, draw=black, fill=white, inner sep=2pt, above of=4] {};
		          \node (6) [circle, draw=black, fill=white, inner sep=2pt, right of=5] {};
			
			     % Draw arrows with labels
		          \draw[->] (1) -- (2) node[midway, below, inner sep=2pt] {$\rho_1^0=0.6$};
		          \draw[->] (2) -- (3) node[midway, left, inner sep=2pt] {$\rho_2^0=0.4$};
		          \draw[->] (2) -- (4) node[midway, below, inner sep=2pt] {$\rho_3^0=0.9$};
                \draw[->] (3) -- (5) node[midway, above, inner sep=2pt] {$\rho_4^0=0.5$};
			    \draw[->] (4) -- (5) node[midway, right, inner sep=2pt] {$\rho_5^0=0.8$};
			    \draw[->] (5) -- (6) node[midway, above, inner sep=2pt] {$\rho_6^0=0.3$};
            \end{scope}
		\end{tikzpicture}
	}
	\end{center}
	\caption{Road network with routing parameter \(\alphaa\) (left) and initial traffic densities (right) of the sample problem}\label{fig: sample network}
\end{figure}
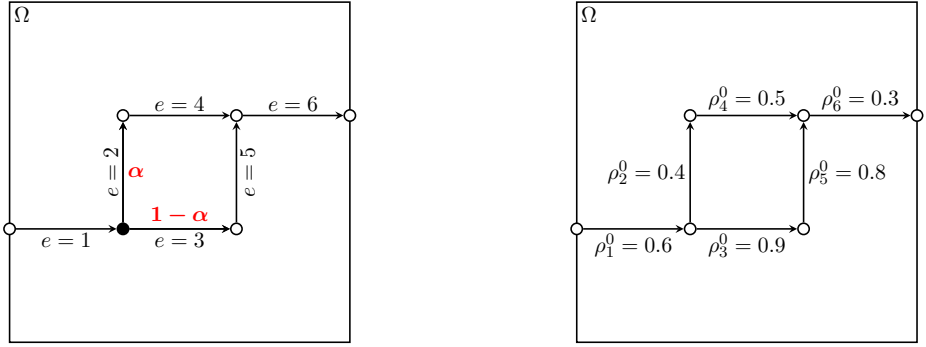

To define the sample problem, we have to fix the parameters of the traffic, emission and dispersion model presented in Section \ref{sec: traffic emission model}.
For comparability, we adopt the same setup as in \cite{Ulke2025}. The corresponding parameter values are summarized in Table~\ref{tab: parameters}. The spatial layout of the road network and the domain \(\Omega\), together with the initial traffic densities, are shown in Figure~\ref{fig: sample network}.
The network contains only one diverging junction, where the road \(e=1\) splits into the roads \(e=2\) and \(e=3\).
Hence, the routing strategy consists of a single routing parameter \(\alpha = \alpha_{2,1}\). Since the routing parameters must sum to one, the second routing parameter is uniquely determined by \(\alpha_{3,1}= 1 - \alpha_{2,1} = 1- \alpha\), see also Figure \ref{fig: sample network} for an illustration.

We consider three control scenarios: optimization of the routing strategy, optimization of the speed limit policy, and joint optimization of both controls.
When only one type of control is optimized, the remaining control variables are fixed. Specifically, we set \(\Vfix = \left(1,1,1,1,1,1\right)^\top\) and \(\alphafix = 0.5\).
%\begin{equation*}
	%\Vfix = \left(1,1,1,1,1,1\right)^\top \quad \text{and} \quad \alphafix = 0.5
%\end{equation*}
The choice of these fixed controls influences the resulting optimal solutions. However, we omit a detailed investigation of this dependence, as it naturally leads to the joint optimization problem.
Finally, we fix the lower and upper bounds of the admissible sets for the optimization to:
\begin{align*}
    \Scal_{\alphaa} & = \Bigl\lbrace \alpha \in \R \mid 0 \le \alpha \le 1 \Bigr\rbrace, \\
    \Scal_{\bm{V}} & = \Bigl\lbrace \Vbm \in \R^{6} \mid 0.25 \leq \Vmax_e \leq 2\Bigr\rbrace, \\
    \Scal & = \Bigl\lbrace (\alpha,\Vbm) \in \R \times \R^6 \mid 
			0 \leq \alpha \leq 1, \,\, 
			0.25 \leq \Vmax_e \leq 2\Bigr\rbrace.
\end{align*}
All single-objective optimization problems are solved using Matlab's \texttt{fmincon}, whereas the multi-objective optimization problems are solved using \texttt{paretosearch}.

\renewcommand{\arraystretch}{1.25}
\begin{table}
    \centering
    \caption{Overview of the parameter values in the traffic, emission, and dispersion model}
    \label{tab: parameters}
    \begin{tabular}{@{}cccccccccccccc@{}}
        \toprule
        \multicolumn{5}{c}{Traffic Model} & & \multicolumn{3}{c}{Emission Model} & &\multicolumn{4}{c}{Dispersion Model} \\
        \cmidrule{1-5} \cmidrule{7-9} \cmidrule{11-14} 
         \(T\) & \(I_e\) & \(\rho_e^{\max}\) & \(q_e^\inn(t)\) & \(\beta_{ji} \) 
            & & \(\Omega\) & \(w_e\) & \(\theta\) 
            & & \(\kappa\) & \(\bm{v}\) & \(\mu\) & \(\phi_0(\x)\) \\
         %\midrule
         \(5\) & \([0,1]\) & \(1\) & \(0.25\) & \(0.5 \) 
            & & \([0,3]^2\) & \(0.1\) & \(0.5\) 
            & & \(0.01\) & \((1,1)^\top\) & \(10^{-6}\) & \(0\)\\
         \bottomrule
    \end{tabular}
\end{table}

\subsection{Maximizing Traffic Efficiency}
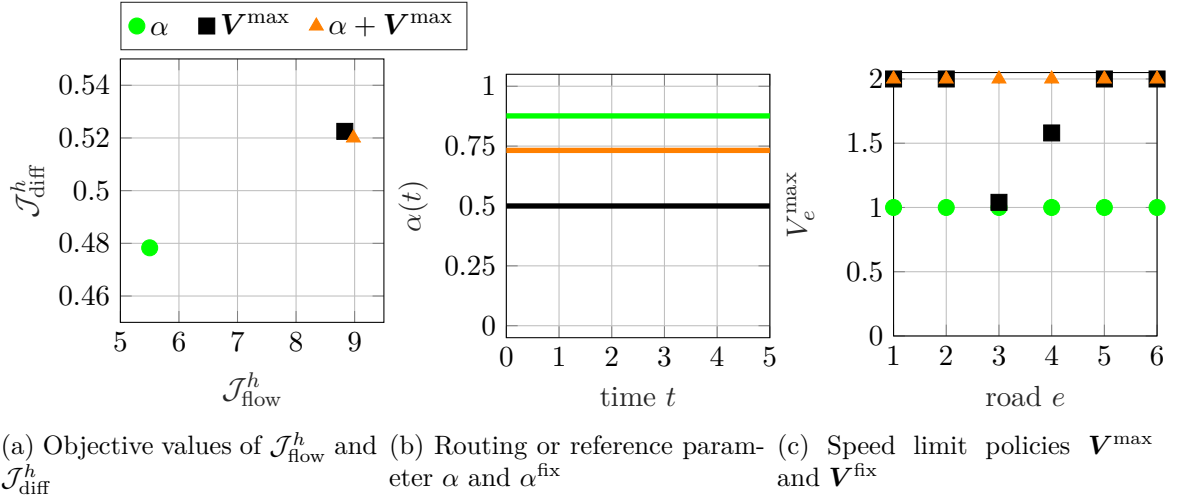
\begin{figure}[h]
	\centering
		\begin{subfigure}[t]{.32\textwidth}
			\centering
\definecolor{mycolor1}{rgb}{1.00000,1.00000,0.00000}%
\begin{tikzpicture}
	
	\begin{axis}[%
		width=0.7\textwidth, % 0.25
		height=0.7\textwidth, % 0.25
		at={(0\textwidth,0\textwidth)},
		scale only axis,
		xmin=5,
		xmax=9.5,
		xtick={5,6,7,8,9},
		xlabel style={font=\color{white!15!black}},
		xlabel={$\J^h_{\flow}$},
		ymin=0.45,
		ymax=0.55,
		ylabel style={font=\color{white!15!black}},
		ylabel={$\J^h_{\diff}$},
		xmajorgrids,
		ymajorgrids,
		axis background/.style={fill=white},
		legend style={
			at={(0.0,1.125)},
			anchor=west,
			legend columns=-1,
			draw=white!15!black,
			/tikz/every even column/.append style={column sep=5pt}
		}
		]
		\addplot [color=green, only marks, mark=*, mark options={solid, fill=green, draw=green}, mark size=3pt]
		table[row sep=crcr]{%
			5.5005	0.4783\\
		};
		\addlegendentry{\(\alpha\phantom{)}\)}
		
%		\addplot [color=red, only marks, mark=triangle*, mark options={solid, fill=red}, mark size=3pt]
%		table[row sep=crcr]{%
%			5.5144	0.4784\\
%		};
%		\addlegendentry{$\alphaa(t)$}
		
		\addplot [color=black, only marks, mark=square*, mark options={solid, fill=black}, mark size=3pt]
		table[row sep=crcr]{%
			8.8323	0.5225 \\
		};
		\addlegendentry{$\Vbm$}
		
		\addplot [color=orange, only marks, mark=triangle*, mark options={solid, fill=orange}, mark size=3pt]
		table[row sep=crcr]{%
			8.9798	0.52\\
		};
		\addlegendentry{\(\alpha+\Vbm\)}
		
%		\addplot [color=blue, only marks, mark=diamond*, mark options={solid, fill=blue}, mark size=3pt]
%		table[row sep=crcr]{%
%			8.9887	0.5191 \\
%		};
%		\addlegendentry{\(\alphaa(t) + \Vbm\)}
		
	\end{axis}
	
	\begin{axis}[%
		width=5.833in,
		height=4.375in,
		at={(0in,0in)},
		scale only axis,
		xmin=0,
		xmax=1,
		ymin=0,
		ymax=1,
		axis line style={draw=none},
		ticks=none,
		axis x line*=bottom,
		axis y line*=left
		]
	\end{axis}
\end{tikzpicture}%
			\caption{Objective values of \( \J_\flow^h\) and \(\J_\diff^h\)}
			\label{fig: flowMax J}
		\end{subfigure}\hfil
		\begin{subfigure}[t]{.32\textwidth}
			\centering
\begin{tikzpicture}

\begin{axis}[%
width=0.7\textwidth, % 0.25
height=0.7\textwidth, % 0.25
at={(0\textwidth,0\textwidth)},
scale only axis,
xmin=0,
xmax=5,
xtick={0,1,2,3,4,5},
xlabel style={font=\color{white!15!black}},
xlabel={time \(t\)},
ymin=-0.05,
ymax=1.05,
ytick={0,0.25,0.5,0.75,1},
ylabel style={font=\color{white!15!black}},
ylabel={\(\alpha(t)\)},
axis background/.style={fill=white},
xmajorgrids,
ymajorgrids,
%legend style={
%	at={(0.01,1.05)},
%	anchor=west,
%	legend columns=-1,
%	draw=white!15!black,
%	/tikz/every even column/.append style={column sep=5pt}
%}
]

\addplot [color=green, line width=2pt]
table[row sep=crcr]{%
	0	0.8761\\
	5	0.8761 \\
};
%\addlegendentry{$\alphaa$}

\addplot [color=black, line width=2pt]
table[row sep=crcr]{%
	0	0.5\\
	5	0.5 \\
};
%\addlegendentry{$\Vbm$}

\addplot [color=orange, line width=2pt]
table[row sep=crcr]{%
	0	0.7321\\
	5	0.7321 \\
};
%\addlegendentry{$\alphaa + \Vbm$}

\end{axis}

\begin{axis}[%
width=5.833in,
height=4.375in,
at={(0in,0in)},
scale only axis,
xmin=0,
xmax=1,
ymin=0,
ymax=1,
axis line style={draw=none},
ticks=none,
axis x line*=bottom,
axis y line*=left
]
\end{axis}
\end{tikzpicture}%
			\caption{Routing or reference parameter \(\alpha\) and \(\alpha^{\mathrm{fix}}\) }
			\label{fig: flowMax routing}
		\end{subfigure}\hfil
		\begin{subfigure}[t]{.32\textwidth}
			\centering
\begin{tikzpicture}
	
	\begin{axis}[%
		width=0.7\textwidth, % 0.25
		height=0.7\textwidth, % 0.25
		at={(0\textwidth,0\textwidth)},
		scale only axis,
		xmin=1,
		xmax=6,
		xtick={1,2,3,4,5,6},
		xlabel style={font=\color{white!15!black}},
		xlabel={road \(e\)},
		ymin=0,
		ymax=2.05,
		ylabel style={font=\color{white!15!black}},
		ylabel={\(\Vmax_e\)},
		axis background/.style={fill=white},
		xmajorgrids,
		ymajorgrids,
		%legend style={
			%	at={(0.01,1.05)},
			%	anchor=west,
			%	legend columns=-1,
			%	draw=white!15!black,
			%	/tikz/every even column/.append style={column sep=5pt}
			%}
		]

		\addplot [color=green, only marks, mark=*, mark options={solid, fill=green, draw=green}, mark size=3pt]
		table[row sep=crcr]{%
			1	1\\
			2	1 \\
			3	1 \\
			4	1 \\
			5	1 \\
			6	1 \\
		};
		%\addlegendentry{$\alphaa$}
		
%		\addplot [color=red, only marks, mark=triangle*, mark options={solid, fill=red}, mark size=2pt]
%		table[row sep=crcr]{%
%			1	1\\
%			2	1 \\
%			3	1 \\
%			4	1 \\
%			5	1 \\
%			6	1 \\
%		};
%		%\addlegendentry{$\alphaa(t)$}
		
		\addplot [color=black, only marks, mark=square*, mark options={solid, fill=black}, mark size=3pt]
		table[row sep=crcr]{%
			1	2\\
			2	2 \\
			3	1.04 \\
			4	1.58 \\
			5	2 \\
			6	2 \\
		};
		%\addlegendentry{$\Vbm$}
		
		\addplot [color=orange, only marks, mark=triangle*, mark options={solid, fill=orange}, mark size=3pt]
		table[row sep=crcr]{%
			1	2\\
			2	2 \\
			3	2 \\
			4	2 \\
			5	2 \\
			6	2 \\
		};
		%\addlegendentry{$\alphaa + \Vbm$}
		
%		\addplot  [color=blue, only marks, mark=diamond*, mark options={solid, fill=blue}, mark size=2pt]
%		table[row sep=crcr]{%
%			1	2\\
%			2	2 \\
%			3	2 \\
%			4	2 \\
%			5	2 \\
%			6	2 \\
%		};
%		%\addlegendentry{$\alphaa(t) + \Vbm$}

	\end{axis}
	
	\begin{axis}[%
		width=5.833in,
		height=4.375in,
		at={(0in,0in)},
		scale only axis,
		xmin=0,
		xmax=1,
		ymin=0,
		ymax=1,
		axis line style={draw=none},
		ticks=none,
		axis x line*=bottom,
		axis y line*=left
		]
	\end{axis}
\end{tikzpicture}%
			\caption{Speed limit policies \(\Vbm\) and \(\Vfix\)}
			\label{fig: flowMax speed}
		\end{subfigure}
		\caption{Influence of routing strategies and speed limits policies on flow maximization}
		\label{fig: flowMax}
	%\end{center}
\end{figure}

First, we investigate the effect of routing strategies and speed limits policies on the efficiency of traffic. Since we measure traffic efficiency by the accumulated traffic flow, this corresponds to maximizing \(\J_\flow^h\). 
The optimal objective values and the corresponding optimal controls are displayed in Figure \ref{fig: flowMax}.
As we also interested in the trade-off between traffic's efficiency and its environmental impact, we also evaluate the mean contamination \(\J_\diff^h\) for the maximizers of \(\J_\flow^h\) and provide the results in Figure \ref{fig: flowMax J}. 
Optimizing only the routing strategy results in the lowest optimal traffic flow with \(\J_\flow^h \approx 5.50\), but also the lowest environmental impact. 
In contrast, optimizing both, the routing strategy and the speed limit policy, yields the highest optimal traffic flow with \(\J_\flow^h \approx 8.98\), whereas optimizing only the speed limit policy results in a slightly lower, yet nearly identical, objective value \(\J_\flow^h \approx 8.83\) compared to optimizing both parameters. 
These results indicate that, for our test case, the speed limit policy has a larger impact on traffic efficiency than the routing strategy.

The optimal routing strategies favor road \(e=2\) over road \(e=3\), see Figure \ref{fig: flowMax routing}.
We can explain this behavior by considering the initial traffic density on the network. The initial density on road \(e=2\) is significantly lower than on road \(e=3\), with values of \(\rho_2^0=0.4\) and \(\rho_3^0=0.9\), respectively. 
Thus, road \(e=2\) has more available capacity for incoming vehicles, whereas directing vehicles towards road \(e=3\) is more likely to induce congestion, which reduces the overall traffic flow.
 
For optimal speed limit policies, we observe a clear tendency towards higher speed limits, see Figure \ref{fig: flowMax speed}, as higher speed limits result in higher traffic flow for the same traffic density.
Optimizing both, routing strategies and speed limit policies, all speed limits are given by their upper bounds. However, when we only optimize the speed limit policy, we observe lower speed limits for the roads \(e=3\) and \(e=4\), while the remaining roads keep the upper bounds.
\vspace{-6cm}

\subsection{Minimizing Traffic Emissions} 
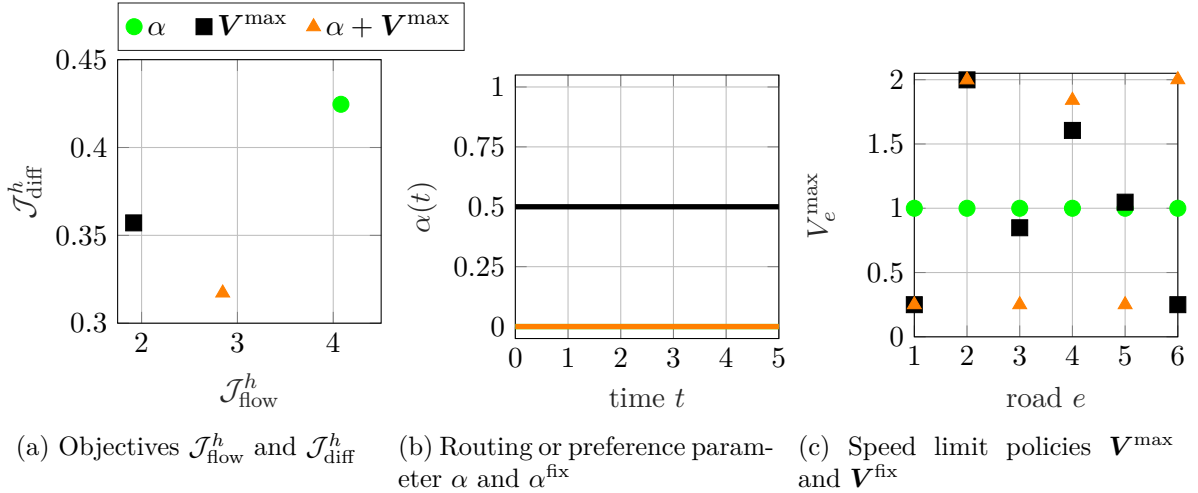
\begin{figure}[h]
	\centering
	\begin{subfigure}[t]{.32\textwidth}
		\centering
\definecolor{mycolor1}{rgb}{1.00000,1.00000,0.00000}%
\begin{tikzpicture}
	
	\begin{axis}[%
		width=0.7\textwidth, % 0.25
		height=0.7\textwidth, % 0.25
		at={(0\textwidth,0\textwidth)},
		scale only axis,
		xmin=1.75,
		xmax=4.5,
		xlabel style={font=\color{white!15!black}},
		xlabel={$\J^h_{\flow}$},
		ymin=0.3,
		ymax=0.45,
		ylabel style={font=\color{white!15!black}},
		ylabel={$\J^h_{\diff}$},
		xmajorgrids,
		ymajorgrids,
		axis background/.style={fill=white},
		legend style={
			at={(0.0,1.125)},
			anchor=west,
			legend columns=-1,
			draw=white!15!black,
			/tikz/every even column/.append style={column sep=5pt}
		}
		]
		
		\addplot [color=green, only marks, mark=*, mark options={solid, fill=green, draw=green}, mark size=3pt]
		table[row sep=crcr]{%
			4.0802	0.4246\\
		};
		\addlegendentry{\(\alpha\phantom{)}\)}
		
%		\addplot [color=red, only marks, mark=triangle*, mark options={solid, fill=red}, mark size=3pt]
%		table[row sep=crcr]{%
%			5.5144	0.4784\\
%		};
%		\addlegendentry{$\alphaa(t)$}
		
		\addplot [color=black, only marks, mark=square*, mark options={solid, fill=black}, mark size=3pt]
		table[row sep=crcr]{%
			1.9174	0.3571 \\
		};
		\addlegendentry{$\Vbm$}
		
		\addplot [color=orange, only marks, mark=triangle*, mark options={solid, fill=orange}, mark size=3pt]
		table[row sep=crcr]{%
			2.8465	0.3172\\
		};
		\addlegendentry{\(\alpha+\Vbm\)}
		
%		\addplot [color=blue, only marks, mark=diamond*, mark options={solid, fill=blue}, mark size=3pt]
%		table[row sep=crcr]{%
%			8.9887	0.5191 \\
%		};
%		\addlegendentry{\(\alphaa(t) + \Vbm\)}
		
	\end{axis}
	
	\begin{axis}[%
		width=5.833in,
		height=4.375in,
		at={(0in,0in)},
		scale only axis,
		xmin=0,
		xmax=1,
		ymin=0,
		ymax=1,
		axis line style={draw=none},
		ticks=none,
		axis x line*=bottom,
		axis y line*=left
		]
	\end{axis}
\end{tikzpicture}%
		\caption{Objectives \( \J_\flow^h\) and \(\J_\diff^h\)}
		\label{fig: pollMin J}
	\end{subfigure}\hfill
	\begin{subfigure}[t]{.32\textwidth}
		\centering
\begin{tikzpicture}
	
	\begin{axis}[%
		width=0.7\textwidth, % 0.25
		height=0.7\textwidth, % 0.25
		at={(0\textwidth,0\textwidth)},
		scale only axis,
		xmin=0,
		xmax=5,
		xtick={0,1,2,3,4,5},
		xlabel style={font=\color{white!15!black}},
		xlabel={time \(t\)},
		ymin=-0.05,
		ymax=1.05,
		ytick={0,0.25,0.5,0.75,1},
		ylabel style={font=\color{white!15!black}},
		ylabel={\(\alpha(t)\)},
		axis background/.style={fill=white},
		xmajorgrids,
		ymajorgrids,
		%legend style={
			%	at={(0.01,1.05)},
			%	anchor=west,
			%	legend columns=-1,
			%	draw=white!15!black,
			%	/tikz/every even column/.append style={column sep=5pt}
			%}
		]

		\addplot [color=green, line width=2pt]
		table[row sep=crcr]{%
			0	0\\
			5	0 \\
		};
		%\addlegendentry{$\alphaa$}
		
		\addplot [color=black, line width=2pt]
		table[row sep=crcr]{%
			0	0.5\\
			5	0.5 \\
		};
		%\addlegendentry{$\Vbm$}
		
		\addplot [color=orange, line width=2pt]
		table[row sep=crcr]{%
			0	0.000018\\
			5	0.000018 \\
		};
		%\addlegendentry{$\alphaa + \Vbm$}

	\end{axis}
	
	\begin{axis}[%
		width=5.833in,
		height=4.375in,
		at={(0in,0in)},
		scale only axis,
		xmin=0,
		xmax=1,
		ymin=0,
		ymax=1,
		axis line style={draw=none},
		ticks=none,
		axis x line*=bottom,
		axis y line*=left
		]
	\end{axis}
\end{tikzpicture}%
		\caption{Routing or preference parameter \(\alpha\) and \(\alpha^{\mathrm{fix}}\)}
		\label{fig: pollMin routing}
	\end{subfigure}\hfill
	\begin{subfigure}[t]{.32\textwidth}
		\centering
\begin{tikzpicture}
	
	\begin{axis}[%
		width=0.7\textwidth, % 0.25
		height=0.7\textwidth, % 0.25
		at={(0\textwidth,0\textwidth)},
		scale only axis,
		xmin=1,
		xmax=6,
		xtick={1,2,3,4,5,6},
		xlabel style={font=\color{white!15!black}},
		xlabel={road \(e\)},
		ymin=0,
		ymax=2.05,
		ylabel style={font=\color{white!15!black}},
		ylabel={\(\Vmax_e\)},
		axis background/.style={fill=white},
		xmajorgrids,
		ymajorgrids,
		%legend style={
			%	at={(0.01,1.05)},
			%	anchor=west,
			%	legend columns=-1,
			%	draw=white!15!black,
			%	/tikz/every even column/.append style={column sep=5pt}
			%}
		]

		\addplot [color=green, only marks, mark=*, mark options={solid, fill=green, draw=green}, mark size=3pt]
		table[row sep=crcr]{%
			1	1\\
			2	1 \\
			3	1 \\
			4	1 \\
			5	1 \\
			6	1 \\
		};
		%\addlegendentry{$\alphaa$}

		\addplot [color=black, only marks, mark=square*, mark options={solid, fill=black}, mark size=3pt]
		table[row sep=crcr]{%
			1	0.25\\
			2	1.9988 \\
			3	0.8483 \\
			4	1.6055 \\
			5	1.0474 \\
			6	0.25 \\
		};
		%\addlegendentry{$\Vbm$}
		
		\addplot [color=orange, only marks, mark=triangle*, mark options={solid, fill=orange}, mark size=3pt]
		table[row sep=crcr]{%
			1	0.25\\
			2	1.9989 \\
			3	0.25 \\
			4	1.8399 \\
			5	0.25 \\
			6	2 \\
		};
		%\addlegendentry{$\alphaa + \Vbm$}

	\end{axis}
	
	\begin{axis}[%
		width=5.833in,
		height=4.375in,
		at={(0in,0in)},
		scale only axis,
		xmin=0,
		xmax=1,
		ymin=0,
		ymax=1,
		axis line style={draw=none},
		ticks=none,
		axis x line*=bottom,
		axis y line*=left
		]
	\end{axis}
\end{tikzpicture}%
		\caption{Speed limit policies \(\Vbm\) and \(\Vfix\)}
		\label{fig: pollMin speed}
	\end{subfigure}
	\caption{Influence of routing strategies and speed limit policies on emission minimization}
	\label{fig: pollMin}
\end{figure}

Next, we analyze the effect of routing strategies and speed limit policies on the environmental impact of traffic, which corresponds to minimizing the mean contamination \(\J_\diff^h\). The results are shown in Figure \ref{fig: pollMin}.

In Figure \ref{fig: pollMin J} we display the optimal objective values when minimizing the mean contamination, while also showing the corresponding values for the accumulated traffic flow.
Again, we observe that optimizing the speed limit policy independently or together with the routing strategy has a greater impact on traffic emissions than routing alone.
In case we only optimize the routing strategy, we obtain the largest optimal value of \(\J_\diff^h \approx 0.42\). When optimizing only the speed limit policy, we obtain a better result with \(\J_\diff^h \approx 0.36\), while optimizing both achieves the best result of \(\J_\diff^h \approx 0.32\). 
Compared to maximizing the traffic flow, the results of optimizing only speed limits and speed limits and routing yields more distinct result. However, as for maximizing traffic efficiency, we conclude for our test case that the speed limit policy has a greater impact on traffic's contribution to air pollution.

Optimizing with the sole objective of lowering the emission of air pollutants creates a strong incentive to keep streets empty for as long as possible, as this achieves zero emission on these unused streets.
This makes the extreme cases \(\alpha = 0\) and \(\alpha = 1\) strong contenders for the optimal choice of the routing strategy. 
Figure \ref{fig: pollMin routing} shows the resulting optimal routing strategies. In both cases, the optimal strategies coincide and are given by \(\alpha =0\).
Notice that when only the speed limit policy is optimized, the routing parameter is prescribed a priori by \(\alpha = 0.5\) and therefore fixed.
The solution \(\alpha = 0\) implies that all vehicles are directed towards road \(e=3\), while road \(e=2\) is closed for new incoming vehicles. At first, this appears counterintuitive since the initial density on road \(e=3\) is higher than on road \(e=2\).
However, we can explain this phenomenon by the structure of the emission rate introduced in Equation \eqref{eq: emission rate}. In our example, the reduction of emissions resulting from lower traffic flows outweighs the increase in emissions caused by higher traffic densities.
Thus, concentrating traffic on the already congested road \(e=3\) leads to lower overall emissions than distributing vehicles in a way that increases flow.
In other words, the network operates in a highly congested regime, with traffic densities already close to road capacity. In this setting, maintaining high densities while reducing flow appears to be more beneficial from an emissions perspective than increasing flow to resolve congestion to reduce densities over time.

Figure \ref{fig: pollMin speed} shows the optimal speed limit polices. 
In the setting where we optimize both, the routing strategy and the speed limit policy, we observe a bang-bang control where speed limits either attain the lower or upper bound.
Further, we notice that roads with lower initial densities are assigned lower speed limits, whereas roads with high initial densities are assigned higher speed limits. Again, as road \(e=2\) is not utilized in this setting, this suggest that it is beneficial to reduce the traffic flow and maintain higher densities than to increase flow in order to resolve congestion. 
When only the speed limit policy is optimized, the pattern changes: now we also have higher speed limits on roads with higher initial density. 
As the routing parameter is fixed to \(\alphafix =0.5\), road \(e=2\) is no longer closed, resulting in a different redistribution of traffic. 
Increasing speed limits on congested roads results in higher traffic flow, which in turn reduces local densities and thus overall emissions.
We also observe that the speed limit on road \(e=1\) is given by the lower bound \(\Vmax_1 = 0.25\), which may be a modeling artifact. A lower speed limit reduces the capacity of the road, limiting inflow and preventing additional congestion from propagating into the network.
\vspace{-6cm}

\subsection{Multi-Objective Optimization}

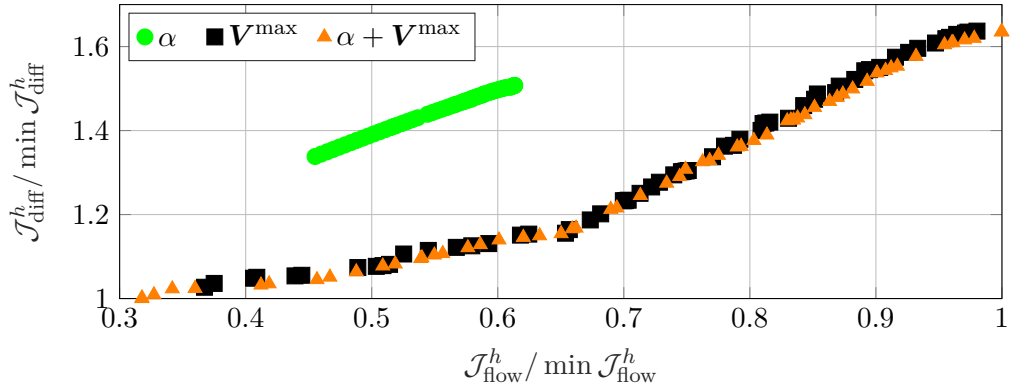
\begin{figure}[h]
	\begin{center}
	\definecolor{mycolor1}{rgb}{1.00000,1.00000,0.00000}%
\begin{tikzpicture}

\begin{axis}[%
width=0.75\textwidth, % 1
height=0.25\textwidth, % 0.788
at={(0\textwidth,0\textwidth)},
scale only axis,
xmin=0.3,
xmax=1,
xlabel style={font=\color{white!15!black}},
xlabel={$\J^h_{\flow}/\min \J^h_{\flow}$},
ymin=1,
ymax=1.7,
ylabel style={font=\color{white!15!black}},
ylabel={$\J^h_{\diff}/\min \J^h_{\diff}$},
xmajorgrids,
ymajorgrids,
axis background/.style={fill=white},
legend style={
	at={(0.01,0.89)},
	anchor=west,
	legend columns=-1,
	draw=white!15!black,
	/tikz/every even column/.append style={column sep=5pt}
}
]
\addplot [color=green, only marks, mark=*, mark options={solid, fill=green, draw=green}, mark size=3pt]
  table[row sep=crcr]{%
0.584739215633425	1.4817504564355\\
0.528069069253577	1.42132099576792\\
0.517334081576582	1.4095246760065\\
0.559490194901944	1.4547564222252\\
0.613433273645642	1.50757458636628\\
0.492740301379705	1.38184482389155\\
0.601916781058205	1.49862333208317\\
0.544295143249808	1.43888559453182\\
0.583615951265023	1.48074474793856\\
0.474149008458308	1.36048810817042\\
0.594359717113799	1.49195540170063\\
0.510798333773417	1.40219350635174\\
0.565774778188874	1.46148589531873\\
0.573755225183589	1.4698653111725\\
0.612525494869896	1.50655442443627\\
0.579781210549349	1.47603226216298\\
0.611496002110599	1.50531532770637\\
0.519503893759833	1.41193072834437\\
0.593804091842694	1.49137776174259\\
0.483544704947534	1.3712822916766\\
0.598363656527022	1.49559655556103\\
0.552015875743257	1.44703392127263\\
0.572925511229152	1.46894300710653\\
0.501834846054721	1.39215346713603\\
0.604981155215777	1.5011773881049\\
0.566723299322024	1.46244744454951\\
0.613163064104336	1.50725827280609\\
0.536395063927272	1.43033184850031\\
0.576380995245778	1.47257498135897\\
0.464667809445121	1.34947233113875\\
0.589628914912657	1.48749903727976\\
0.532223845427687	1.4258815758507\\
0.587419560491991	1.48439392599419\\
0.459850151324146	1.34388105679184\\
0.588078298491104	1.48514622378609\\
0.497338608779076	1.38703180037919\\
0.603655398572352	1.4999452993699\\
0.563055921715633	1.45858280843872\\
0.613340720578566	1.50746261892058\\
0.478853973215834	1.36591507558559\\
0.59617753901122	1.49391026522762\\
0.548256158223713	1.4430121512963\\
0.57933425824513	1.47547580142527\\
0.577286117513509	1.4734002537333\\
0.612094860438182	1.50601698456373\\
0.515170478237765	1.40709907512741\\
0.58933473924701	1.48700388806828\\
0.523857731009677	1.41667264907801\\
0.455038572122824	1.3382358597504\\
0.580788047051033	1.47794258852006\\
0.610783307850967	1.50436793977033\\
0.555740001065209	1.45095797980695\\
0.609854557458907	1.50295345996858\\
0.488146828159643	1.37659217423586\\
0.600217281785097	1.49720362946386\\
0.570403441881186	1.4661872072141\\
0.612925206260161	1.5069614679173\\
0.506308446920305	1.39721113519429\\
0.606369949839264	1.50223049042511\\
0.469445368808664	1.35500684172393\\
};
\addlegendentry{\(\alpha\phantom{)}\)}

\addplot [color=black, only marks, mark=square*, mark options={solid, fill=black}, mark size=3pt]
  table[row sep=crcr]{%
0.971129605264144	1.63220599356974\\
0.439122432381381	1.05397382937397\\
0.980177722896926	1.63681721525132\\
0.367340870984157	1.02691970502691\\
0.964094553696386	1.62901039767838\\
0.973013333681832	1.63445944935736\\
0.88323114258764	1.52174917031704\\
0.808876050917289	1.40030747128441\\
0.503650278602628	1.07673016274595\\
0.745704121550676	1.30201100193502\\
0.69961227912338	1.23291803097704\\
0.722066395746212	1.26572437379969\\
0.67339570801898	1.18719947582886\\
0.779913690437725	1.36287888470694\\
0.579393150720382	1.12556799091613\\
0.842702248894237	1.4593893440256\\
0.891157019805789	1.54233408960202\\
0.792097190458414	1.37928433489507\\
0.770378056798388	1.33762038354799\\
0.955486903766884	1.61871274516515\\
0.94729083874014	1.60843574567181\\
0.903012443826861	1.55044061872819\\
0.44452854628606	1.05538646690897\\
0.830668400529112	1.4293120753689\\
0.853586742054213	1.48714790458162\\
0.514154287106645	1.08140437962063\\
0.810468878136924	1.41662646080139\\
0.868051434455338	1.49024870745582\\
0.749365399794845	1.3037555638872\\
0.617992506646953	1.15107509093033\\
0.787010399993788	1.36483490779058\\
0.712874584374063	1.25053832873565\\
0.488992436077072	1.07376345795364\\
0.592822358815812	1.130860487052\\
0.728240726923029	1.27701093914663\\
0.408827499344421	1.05060679204521\\
0.544879148428419	1.11461976943064\\
0.406242607798659	1.04842115409085\\
0.85136069704678	1.47458249121574\\
0.815844149705118	1.42012144592834\\
0.656762938069576	1.16413487508818\\
0.915531112744959	1.57466281200701\\
0.701217464523559	1.23312002132857\\
0.751126847399322	1.30488067869886\\
0.653717688636675	1.15560780629662\\
0.525360224453046	1.10606599316589\\
0.624600534044552	1.15343053815477\\
0.681631467943521	1.20167448139604\\
0.958433897680886	1.62198191100249\\
0.375000756784884	1.03602247684637\\
0.925970198707418	1.58620939112316\\
0.870902756138536	1.50582655544207\\
0.739617060832145	1.29451622428625\\
0.508976793277097	1.07839483107853\\
0.7034527053053	1.23420548400297\\
0.74708813366909	1.30264614879804\\
0.894955111496908	1.54513447747545\\
0.567116378667766	1.12191156033663\\
0.933202550283104	1.59588926548918\\
0.812633313380043	1.41801736352978\\
};
\addlegendentry{$\Vbm$}

\addplot [color=orange, only marks, mark=triangle*, mark options={solid, fill=orange}, mark size=3pt]
  table[row sep=crcr]{%
0.999473771241034	1.63497542004316\\
1	1.63532414029511\\
0.327205876630045	1.00906231831348\\
0.713124568108673	1.24541971645648\\
0.851289365526977	1.45531166922343\\
0.317741104424618	1.00153103726677\\
0.418612039584411	1.03513019098909\\
0.694456195450064	1.21567618550726\\
0.931824899890352	1.57659021996918\\
0.829659292811091	1.42319118074408\\
0.456592450034205	1.04509735755495\\
0.863251229391451	1.46938947065644\\
0.790002319761367	1.36101756287859\\
0.813558117999912	1.38944696324839\\
0.31762901382601	1\\
0.539844564920388	1.09702119616889\\
0.508762462312616	1.07759017491991\\
0.762530163691497	1.32579643247722\\
0.917059333054659	1.55335967672713\\
0.912072324575805	1.54933882205761\\
0.74433451041775	1.29050934380402\\
0.837875894953349	1.43011413368078\\
0.775112155216941	1.34109553912394\\
0.733995331266549	1.27472717478084\\
0.576287428081689	1.12113713889909\\
0.869255025074942	1.47893357506466\\
0.970371166751539	1.61592660143158\\
0.960353389320675	1.60928204810439\\
0.881568393712112	1.49916276305396\\
0.659157712335887	1.16621519142592\\
0.893055117059492	1.51733223614209\\
0.873776522182457	1.48691129145949\\
0.359825952654997	1.02412275563884\\
0.586171318650909	1.1284748399182\\
0.834642211388577	1.42510476399387\\
0.412147861538358	1.03222850098418\\
0.600973126441563	1.13959525227573\\
0.803196881693246	1.37650749027028\\
0.538919414303176	1.09533497589565\\
0.793101572360517	1.36236706199481\\
0.466817220086707	1.05109084121235\\
0.549706424109008	1.10378441747396\\
0.768125746093863	1.32859056256975\\
0.518662402808048	1.08233893169626\\
0.633215005265015	1.14949780486526\\
0.341774618471564	1.02321482291037\\
0.900900358936803	1.5372076225195\\
0.6621443878748	1.16782548482749\\
0.556356330980981	1.10718093398772\\
0.907646853327373	1.5432716693138\\
0.978016384618706	1.61969963769482\\
0.487754423783147	1.06380446072537\\
0.620269458351056	1.14569919485426\\
0.843603686864182	1.43787291804327\\
0.689652701347139	1.2125144619373\\
0.954435402802264	1.60428463031239\\
0.832532411960234	1.42503943986181\\
0.961144498753409	1.6094848698547\\
0.650292852502688	1.15454916074426\\
0.749020307134366	1.30788795839693\\
};
\addlegendentry{\(\alpha+\Vbm\)}

\end{axis}

\begin{axis}[%
width=5.833in,
height=4.375in,
at={(0in,0in)},
scale only axis,
xmin=0,
xmax=1,
ymin=0,
ymax=1,
axis line style={draw=none},
ticks=none,
axis x line*=bottom,
axis y line*=left
]
\end{axis}
\end{tikzpicture}%
	\caption{Influence of routing strategies and speed limit policies on the trade-off between emission reduction and flow maximization}
	\label{fig: Pareto Front}
	\end{center}
\end{figure}

Having studied the individual objectives separately, we now investigate the trade-off between traffic efficiency and environmental impact through multi-objective optimization. The results of the single-objective optimization problems indicate that these objectives are conflicting, which means increasing traffic efficiency leads to higher emissions, while reducing emissions comes at the expense of lower traffic efficiency, cf. Figures~\ref{fig: flowMax J} and \ref{fig: pollMin J}.
To analyze this trade-off, we employ multi-objective optimization and determine Pareto-optimal solutions. Again, we consider optimizing the routing strategy, the speed limit policy, and both controls simultaneously. The resulting Pareto fronts are shown in Figure~\ref{fig: Pareto Front}. For better comparability, the axes are normalized by the optimal values of \(\J_{\diff}^h\) and \(\J_{\flow}^h\), respectively.

The results confirm the observations from our previous analysis: the speed limit policy has a substantially larger impact on the optimal solution than the routing strategy.
Additionally, it produces the better solution as the Pareto front obtained by optimizing speed limits dominates the front obtained by optimizing routing. 
 Moreover, the Pareto front corresponding to the joint optimization of the routing strategy and speed limit policy is very close to that obtained by optimizing only the speed limit policy, indicating that routing provides only a limited additional benefit in the considered example.
The intuition behind this observation is that speed limits directly affect traffic dynamics throughout the entire network, whereas routing decisions are made only at a single diverging junction. Although routing influences the global traffic dynamic, its effect is inherently more localized than that of the speed limit policy.

\vspace{-6cm}
\section{Conclusion and Outlook}

We have investigated the influence of routing strategies and speed limit policies on optimal solutions in a macroscopic traffic emission model. Numerical experiments on a sample road network yield two main findings:

First, speed limits dominate the routing. In both the single-objective and multi-objective settings, optimizing speed limits alone achieves a larger reduction in mean contamination and a larger gain in accumulated traffic flow than optimizing only the routing strategy. 
This is consistent with the global nature of speed limit policies: a speed limit reduction on a single road simultaneously lowers the flux capacity of that road everywhere and at all times, while a routing adjustment has only a local effect that propagates through the network over the given time horizon.

Second, optimizing the routing strategy and the speed limit policy simultaneously yields the best objective values in all scenarios, but the additional gain over controlling the speed limit policy alone is moderate. This suggests that speed limits are the primary lever, and routing provides fine-tuning.

\medbreak
In future work, it would be of interest to consider second-order traffic models such as \cite{AwRascle2000,Zhang2002}, since the LWR model used here is not able to capture stop-and-go waves.
Additionally, second-order traffic models allow to incorporate information about the velocity and acceleration of vehicles in the emission model and was already explored in \cite{Balzotti2022}. 
Furthermore, allowing time-dependent routing strategies would better reflect real-world traffic management, where routing recommendations can be adapted dynamically in response to changing traffic conditions.
However, to avoid unrealistic instantaneous switching, minimum dwell-time constraints should be incorporated.
Finally, extending the analysis to larger and more realistic networks is an important direction for future research, as city-scale networks with real traffic data would test the robustness of the findings.

%\begin{small}
%    \noindent  \textbf{Use of AI tools declaration}
%    The authors used Claude (Anthropic) and ChatGPT (OpenAI) as a language tool to check grammar and improve the clarity of the manuscript text. No content was generated by AI; all scientific ideas, results, and interpretations are solely the authors' own work.
%\end{small}
%\begin{acknowledgement}
%  An acknowledgement may be placed at the end of the article.
%\end{acknowledgement}

\vspace{\baselineskip}
%% The style of the following references should be used in all documents.


\begin{thebibliography}{1}
\bibitem{Alvarez2017}% 
	L. J. Alvarez-Vázquez, N. García-Chan, A. Martínez, and M. E. Vázquez-Méndez, \textit{Numerical simulation of air pollution due to traffic flow in urban networks}, J Comput Appl Math \textbf{326}, pp. 44--61 (2017).

\bibitem{Alvarez2018}% 
	L. J. Alvarez-Vázquez, N. García-Chan, A. Martínez, and M. E. Vázquez-Méndez, \textit{Optimal control of urban air pollution related to traffic flow in road networks}, Math Control Relat Fields \textbf{8}, pp. 177--193 (2018). 

\bibitem{AwRascle2000}% 
	A. Aw, and M. Rascle, \textit{Resurrection of “second order” models of traffic flow}, SIAM J. Appl. Math. \textbf{60}, pp. 916--938 (2000).

\bibitem{Balzotti2022}% 
	B. Balzotti, M. Briani, B. De Filippo, and B. Piccoli, \textit{A computational modular approach to evaluate \(\mathrm{NO}_x\) emissions and ozone production due to vehicular}, Discrete Contin Dyn Syst B \textbf{27}, pp. 3455--3486 (2022).

\bibitem{Berrone2012}% 
	S. Berrone, F. De Santi, S. Pieraccini, and M. Marro, \textit{Coupling traffic models on networks and urban dispersion models for simulating sustainable mobility strategies}, Comput Math Appl \textbf{64}, pp. 1975--1991 (2012).
    
\bibitem{Branke2008}% 
	J. Branke, D. Kalyanmoy, K. Miettinen, and R. Slowinski, \textit{Mulitobjective optimization}, Springer, Berlin (2008).

\bibitem{Bressan2014}% 
	A. Bressan, S. Čanić, M. Garavello, M. Herty, and B. Piccoli, \textit{Flows on networks: recent results and perspectives}, EMS Surv. Math. Sci. \textbf{1}, pp. 47--111 (2014).
   

\bibitem{Ehrgott2005}% 
	M. Ehrgott, \textit{Multicriteria optimization}, Springer, Berlin (2005).

\bibitem{Piccoli2006}% 
	M. Garavello, and B. Piccoli, \textit{Traffic flow on networks}, American Institute of Mathematical Sciences, Springfield (2006).

\bibitem{Garavello2016}% 
	M. Garavello, K, Han, and B. Piccoli, \textit{Models for vehicular traffic on networks}, Am Inst Math Sci (2016).

\bibitem{Ulke2025}% 
	S. G\"ottlich, M. Herty, and A. Ulke, \textit{Speed limits in traffic emission models using multi-objective optimization}, Optim. Eng. \textbf{26}, pp. 199--227 (2025).

\bibitem{Goettlich2015}% 
	S. G\"ottlich, M. Herty, and U. Ziegler,
    \textit{Modeling and optimizing traffic light settings in road networks},
 	Comput Oper Res \textbf{55}, pp. 36--51 (2015).

\bibitem{Goatin2016}% 
	P. Goatin, S. G\"ottlich, and O. Kolb,
    \textit{Speed limit and ramp meter control for traffic flow networks},
 	Eng Optim \textbf{48}, pp. 1121--1144 (2016).

\bibitem{Goatin2026}% 
	P. Goatin, A. Klar, C. Mezquita-Nieto,
    \textit{Discrete adjoint gradient computation for multiclass traffic flow models on road networks},
 	arXiv preprint, (2026).

\bibitem{Greenshields1935}% 
	B. D. Greenshields, J. R. Bibbins, W. S. Channing, and H. H. Miller, \textit{A study of traffic capacity}, Highway research board proceedings, National Research Council (USA), Highway Research Boa (1935).

\bibitem{Herty2009}% 
	M. Herty, J. P. Lebacque, and S. Moutari, \textit{A novel model for intersections of vehicular traffic flow}, Netw Heterog Media \textbf{4}, pp. 813--826 (2009).
    
\bibitem{LighthilWhitham1955}% 
	M. J. Lighthill, and G. Whitham, \textit{On kinematic waves. II. A theory of traffic flow on long crowded roads}, Proc R Soc Lond Ser A Math Phys Sci \textbf{229}, pp. 317--345 (1955).

\bibitem{Richards1956}% 
	P. I. Richards, \textit{Shock waves on the highway}, Oper. Res. \textbf{4}, pp. 42--51 (1956).

\bibitem{Skiba2000}% 
	Y. N. Skiba, and D. Parra-Guevara, \textit{Industrial pollution transport. Part 1. Formulation of the problem and air pollution estimates}, Environ Model Assess \textbf{5}, pp.169--175 (2000).

\bibitem{Stockie2011}% 
	J. M. Stockie, \textit{The mathematics of atmospheric dispersion modeling}, SIAM Rev \textbf{53}, pp. 349--372 (2011).

\bibitem{Treiber2013}% 
	M. Treiber, and A. Kesting, \textit{Traffic flow dynamics: data, models, and simulation}, Springer, Berlin (2013).

\bibitem{Alvarez2019}% 
	M. E. Vázquez-Méndez, L. J. Alvarez-Vázquez, N. García-Chan, and A. Martínez, \textit{Optimal management of an urban road network with an environmental perspective}, Comput Math Appl \textbf{77}, pp. 1786--1797 (2019).


\bibitem{WHO2021}% 
	World Health Organization, \textit{WHO global air quality guidelines: particulate matter (PM2.5 and PM10), ozone, nitrogen dioxide, sulfur dioxide and carbon monoxide}, WHO, Geneva (2021).

\bibitem{Zhang2002}% 
	H. M. Zhang, \textit{A non-equilibrium traffic model devoid of gas-like behavior}, Transp. Res. B \textbf{36}, pp. 275--290 (2002).



\end{thebibliography}
\end{document}